\documentclass[12pt,a4paper,leqno]{article}
\usepackage[latin2]{inputenc}
\usepackage{amsmath,amsthm}
\usepackage{amsfonts}
\usepackage{amssymb}
\usepackage{graphicx}
\usepackage[left=2.50cm, right=2.50cm, top=2.50cm, bottom=2.50cm]{geometry}
\usepackage{sectsty}
\usepackage{titlesec}

\titlelabel{\thetitle.\quad}
\sectionfont{\centering}
\subsectionfont{\centering}

\DeclareMathOperator{\ind}{ind}
\DeclareMathOperator{\res}{res}

\def\Q{{\mathbb Q}}
\def\N{{\mathbb N}}
\def\Z{{\mathbb Z}}
\def\R{{\mathbb R}}
\def\C{{\mathbb C}}

\newtheoremstyle{AMH}{3pt}{3pt}{\itshape}{\parindent}{\scshape}{.}{.5em}{}

\theoremstyle{AMH}
\newtheorem{lemma}{Lemma}
\newtheorem{theorem}[lemma]{Theorem}
\newtheorem{corollary}[lemma]{Corollary}

\newenvironment{prf}[1][\proofname]{
  \proof[\scshape #1]
}{\endproof}

\title{
\textbf{A GENERALIZATION OF SIMPLEST NUMBER FIELDS AND THEIR INTEGRAL BASIS}
}
\author{L. REMETE
\thanks{
        Supported through the \'UNKP-19-3 New National Excellence Program of the Ministry for Innovation and Technology.} \\
{\small University of Debrecen, Mathematical Institute} \\
{\small H--4002 Debrecen Pf.400., Hungary,}\\ 
{\small e--mail: remete.laszlo@science.unideb.hu}
}
\date{}

\renewenvironment{abstract}
 {\par\noindent\textbf{\abstractname.}\ \ignorespaces}
 {\par\medskip}

\begin{document}
\maketitle
\thispagestyle{empty}

\renewcommand{\thefootnote}{\arabic{footnote}}
\setcounter{footnote}{0}
\begin{abstract}
An integral basis of the simplest number fields of degree 3,4 and 6 over $\Q$ are well-known, and widely investigated. We generalize the simplest number fields to any degree, and show that an integral basis of these fields is repeating periodically.
\end{abstract}

\section{Introduction}
Let $a,b,c,d\in\Q$ and 
$$\sigma:\C\mapsto\C, \qquad \sigma:z \mapsto \frac{az+b}{cz+d}.$$
Assume, that $f(X)\in\Z[X]$ is a polynomial with real roots, such that $\sigma$ transitively permutes the roots of $f(X)$. In this case, if $\beta$ is a root of $f(X)$, then $\Q(\beta)$ is a totally real cyclic number field. 
This requires, that the matrix 
$$M=\left(
\begin{array}{rr}
a & b\\
c & d
\end{array}
\right)\in PGL_2(\Q),$$
is of finite order. It can be shown that each non-trivial torsion element of $PGL_2(\Q)$ has order 2,3,4 or 6, whence, there exist such polynomials and number fields only of degrees 2,3,4 and 6.\\
According to  D.Shanks \cite{shanks}, A.J.Lazarus \cite{laz}, G.Lettl, A.Pethő and P.Voutier \cite{lpv} and A.Hoshi  \cite{hoshi6} the polynomials of degrees 3,4 and 6 with these properties are called \textit{simplest} polynomials, and the corresponding number fields are called \textit{simplest} number fields. These fields have an extensive literature. \\
First, the simplest cubic fields are investigated by H.Cohn \cite{cohn} and D.Shanks \cite{shanks}, because they have easily computable and relatively large class numbers.\\
M.N.Gras \cite{gras1}, \cite{gras2}, V.Ennola \cite{ennola1}, \cite{ennola2} and A.J.Lazarus \cite{laz} investigated the unit group of the simplest fields. K.Foster \cite{ht90} obtained the simplest parametric polynomials, just by using a special identity of units in cyclic extensions. It shows, that these fields have some other unique and interesting properties.\\
E.Thomas \cite{thomas}, M.Mignotte \cite{mig}, G.Lettl, A.Pethő and P.Voutier \cite{lpv}, \cite{lpv2} and G.Lettl and A.Pethő \cite{lp}, I.Gaál \cite{book} solved Thue equations corresponding to the simplest polynomials in absolute case, and C.Heuberger \cite{heu}, I.Gaál, B.Jadrijevi\'{c} and L.Remete \cite{gjr} in certain relative cases. \\
A.Hoshi  \cite{hoshi3}, \cite{hoshi4}, \cite{hoshi6} gave a correspondence between solutions of a family of Thue equations and the isomorphism classes of the simplest number fields. He extended his results to a family of polynomials of degree 12, which has similar properties as the simplest polynomials, but over $\Q(\sqrt{-3})$. This generalization provided the main motivation of our approach.\\
Assume that $\beta$ is a root of a simplest polynomial $f(X)\in\Z[X]$ of degree $n$, and the Möbius transformation $\sigma$ permutes its roots transitively.\\
Then the conjugates of $\beta$ are:
$$\{ \beta, \sigma(\beta), \sigma^2(\beta), \ldots , \sigma^{n-1}(\beta) \}.$$
Let $m$ denote
$$m=\frac{\beta+\sigma(\beta)+\sigma^2(\beta)+\ldots+\sigma^{n-1}(\beta) }{n}.$$
If $n=3$ or $n=6$, and 
$$\sigma:\C\mapsto\C, \qquad \sigma:z \mapsto \frac{az-1}{z+a+1}.$$
where $a=0$ or $-1$ for degree 3, and $a=1$ or $-2$ for degree 6, then  we can write the minimal polynomial of $\beta$ over $\Q$ in the form:
$$f^{(3)}_m(X)=X^3-3mX^2-3(m+1)X-1,$$
$$f^{(6)}_m(X)=X^6-6mX^5-15(m+1)X^4-20X^3+15mX^2+6(m+1)X+1.$$
Therefore, if $m$ is a rational parameter, then the polynomials above are cyclic, they have rational coefficients and 
$$\sigma:z \mapsto \frac{az-1}{z+a+1}$$
permutes their roots transitively.\\
Although these polynomials have nice properties for any $m\in\Q$, we say that their root generates a simplest number field, if they are irreducible, and have integer coefficients.\\
In Section \ref{gen} we will generalize these polynomials to any degree. We keep the property that the coefficients are rational numbers, and a Möbius transformation
$$\sigma:z \mapsto \frac{\alpha z-1}{z+\alpha+1}$$
permutes the roots transitively, where $\alpha\in\C$ is an algebraic integer, such that 
$$M=\left(
\begin{array}{rr}
\alpha & -1\\
1 & \alpha+1
\end{array}
\right)\in PGL_2(\C),$$
has finite order. In such way, we will obtain infinite parametric families of polynomials of degree $n$, such that their Galois groups are cyclic over $\Q(\alpha)$. \\
We get some well-known families as special cases.
\begin{itemize}
\item Simplest cubic fields: $n=3$, $m=\frac{t}{3}$, $t\in \Z$ (c.f. \cite{shanks})
\item Simplest sextic fields: $n=6$, $m=\frac{t}{3}$, $t\in \Z\setminus\{-8,-3,0,5\}$ (c.f. \cite{gras2},\cite{sextic})
\item Family of polynomials of degree 12 investigated by A.Hoshi: $n=12$, $m=\frac{t}{3}$, \\
$t\in \Z\setminus\{-8,-3,0,5\}$ (c.f. \cite{hoshi12})
\end{itemize}
We call the number fields generated by a root of our polynomials \textit{generalized simplest} number fields. As the above examples show, these are indeed generalizations of the simplest fields of degrees 3 and 6, but we remark, that for degree 4, our family of number fields is different from the simplest quartic fields.\\
In Section \ref{Sdisc} and \ref{Sired}, we describe the discriminants of these polynomials, and give a sufficient condition for the irreducibility.\\\\
Explicit integral bases of the simplest fields are given by D.Shanks \cite{shanks}, H.K.Kim and J.H.Lee \cite{kl}, G.Lettl, A.Pethő, P.Voutier \cite{lpv2} and I.Gaál and L.Remete \cite{sextic}. In each cases an integral basis of these fields is repeating periodically in $m$.
In Section \ref{prib}, we will show, that in our generalization the same phenomena occurs. Let $t$ be an integer. If $t$ is congruent to 1 or 2 modulo 3, then we set $m=t$, otherwise we set $m=\frac{t}{3}$. In such way we obtain families of parametric polynomials, which have integer coefficients for any integer parameter $t$. We prove, that if the $3$-free part of $t^2+t+1$ (or $t^2+3t+9$, respectively) is square-free, then an integral basis of the fields generated by a root $\beta_t$ of $f^{(n)}_t(X)$ is repeating periodically.\\
We say that an integral basis is repeating periodically modulo $n_0$, if for $r=0,\ldots,n_0-1$ there are polynomials $h^{(r)}_i(X)\in\Q[X]$, $(i=0,\ldots n-1)$, such that, if $t\equiv r\mod{n_0}$ and $f^{(n)}_t(X)$ is irreducible, then  
$$\left(h^{(r)}_0(\beta_t),h^{(r)}_1(\beta_t),h^{(r)}_2(\beta_t),\ldots,h^{(r)}_{n-1}(\beta_t)\right), $$
is an integral basis of $\Q(\beta_t)$. (See \cite{ibm},\cite{sextic},\cite{ibpf}).\\
We give a general upper bound for the period length for any $n\in\N$, and for $n\leq12$ we significantly improve it. Using these better bounds, we are able to determine the smallest period lengths for $n=2,3,4,5,6,8,9,12$.

\section{Generalized simplest number fields}\label{gen}
In this section we define the parametric family $f^{(n)}_m$ of polynomials of degrees $n>1$, such that there exists a Möbius transformation
$$\sigma:z \mapsto \frac{\alpha z-1}{z+\alpha+1},\qquad \alpha\in\C$$
which permutes the roots of $f^{(n)}_m(X)$ transitively for any parameter $m$.
All along this section let $m$ be a rational parameter.\\
Let $g,h:\N\mapsto\Q$ be defined by\\
$$g(i):=\left\lbrace\begin{array}{lr}
1, & \mbox{ if } i\equiv 0 \mod{6},\\
-m, & \mbox{ if } i\equiv 1 \mod{6},\\
-m-1, & \mbox{ if } i\equiv 2 \mod{6},\\
-1, & \mbox{ if } i\equiv 3 \mod{6},\\
m, & \mbox{ if } i\equiv 4 \mod{6},\\
m+1, & \mbox{ if } i\equiv 5 \mod{6}.
\end{array}\right.\qquad\mbox{and}\qquad
h(i):=\left\lbrace\begin{array}{lr}
0, & \mbox{ if } i\equiv 0 \mod{6},\\
-1, & \mbox{ if } i\equiv 1 \mod{6},\\
-1, & \mbox{ if } i\equiv 2 \mod{6},\\
0, & \mbox{ if } i\equiv 3 \mod{6},\\
1, & \mbox{ if } i\equiv 4 \mod{6},\\
1, & \mbox{ if } i\equiv 5 \mod{6}.
\end{array}\right.$$
For $n\geq 0$, let
$$f^{(n)}_m(X):=\sum_{i=0}^{n}\binom{n}{i}\cdot X^i \cdot g(n-i)\qquad\in\Q[m,X],$$
and
$$r^{(n)}(X):=\sum_{i=0}^{n}\binom{n}{i}\cdot X^i \cdot h(n-i)\qquad\in\Z[X].$$
These polynomials satisfy the following recursion formulas:
\begin{lemma}\label{fr}
For $n\geq0$,
\begin{equation}\label{fn+1}
f^{(n+1)}_m(X)=(X-m)\cdot f^{(n)}_m(X)+(m^2+m+1)\cdot r^{(n)}(X),
\end{equation}
and
\begin{equation}\label{rn+1}
r^{(n+1)}(X)=(X+m+1)\cdot r^{(n)}(X)-f^{(n)}_m(X).
\end{equation}

\end{lemma}

\begin{prf}
We prove only (\ref{rn+1}), the other identity can be proved in the same way. We compare the coefficients of $X^k$ in the two side of the equation.\\
The coefficient of $X^k$ on the left hand side of (\ref{rn+1}) is 
$$\binom{n+1}{k}\cdot h(n+1-k)$$
and on the right side of (\ref{rn+1}) is 
$$\binom{n}{k-1}\cdot h(n+1-k)+(m+1)\cdot \binom{n}{k}\cdot h(n-k)-\binom{n}{k}\cdot g(n-k).$$
It is easy to check that for any remainder of $n-k \mod{6}$, these coefficients are the same.
\end{prf}
By the definition of $f^{(n)}_m$ and $r^{(n)}$, it is easy to determine their derivatives.
\begin{lemma}\label{deri}
For $n\geq0$,
$$\left(f^{(n+1)}_m\right)'(X)=\sum_{i=1}^{n+1}\binom{n+1}{i}\cdot i\cdot X^{i-1} \cdot g(n+1-i)=(n+1)\sum_{i=1}^{n+1}\binom{n}{i-1}\cdot X^{i-1} \cdot g(n+1-i)=$$
$$=(n+1)\cdot f^{(n)}_m(X),$$
and similarly
$$\left(r^{(n+1)}\right)'(X)=(n+1)\cdot r^{(n)}(X).$$
\end{lemma}
Now we prove a nice variable transformation identity for $f^{(n)}_m(X)$.

\begin{theorem}\label{ff}
For any $\alpha\in\C$ and $n\geq1$,
\begin{equation}\label{vtrans}
(X+\alpha+1)^n\cdot f^{(n)}_m\left(\frac{\alpha X-1}{X+\alpha+1}\right)=f^{(n)}_m(\alpha)\cdot f^{(n)}_m(X)-(m^2+m+1)\cdot r^{(n)}(\alpha)\cdot r^{(n)}(X).
\end{equation}
\end{theorem}

\begin{prf}
If $n=1$, then $f^{(1)}_m(X)=X-m$, $r^{(1)}(X)=-1$, and
$$(X+\alpha+1)\cdot \left(\frac{\alpha X-1}{X+\alpha+1}-m\right)=(\alpha-m)\cdot(X-m)-(m^2+m+1)\cdot(-1)\cdot(-1).$$
Hence, (\ref{vtrans}) is true for $n=1$. \\
Assume, that (\ref{vtrans}) is true for $n$:
\begin{equation}\label{base}
(X+\alpha+1)^n\cdot f^{(n)}_m\left(\frac{\alpha X-1}{X+\alpha+1}\right)=f^{(n)}_m(\alpha)\cdot f^{(n)}_m(X)-(m^2+m+1)\cdot r^{(n)}(\alpha)\cdot r^{(n)}(X).
\end{equation}
Let
$$F_{n+1}(X):=\frac{f^{(n+1)}_m(\alpha)\cdot f^{(n+1)}_m(X)-(m^2+m+1)\cdot r^{(n+1)}(\alpha)\cdot r^{(n+1)}(X)}{(X+\alpha+1)^{n+1}}.$$
Using Lemma \ref{deri}, the derivative of $F_{n+1}$ with respect to $X$ is:
$$F'_{n+1}(X)=\frac{n+1}{(X+\alpha+1)^{n+2}}\cdot \left[\left(f^{(n+1)}_m(\alpha)\cdot f^{(n)}_m(X)-(m^2+m+1)\cdot r^{(n+1)}(\alpha)\cdot r^{(n)}(X)\right)\cdot(X+\alpha+1)-\right.$$
$$\left.f^{(n+1)}_m(\alpha)\cdot f^{(n+1)}_m(X)+(m^2+m+1)\cdot r^{(n+1)}(\alpha)\cdot r^{(n+1)}(X)\right].$$
By Lemma \ref{fr}, substitute
$$r^{(n+1)}(X)=(X+m+1)\cdot r^{(n)}(X)-f^{(n)}_m(X),$$
$$f^{(n+1)}_m(X)=(X-m)\cdot f^{(n)}_m(X)+(m^2+m+1)\cdot r^{(n)}(X),$$
$$r^{(n+1)}(\alpha)=(\alpha+m+1)\cdot r^{(n)}(\alpha)-f^{(n)}_m(\alpha),$$
$$f^{(n+1)}_m(\alpha)=(\alpha-m)\cdot f^{(n)}_m(\alpha)+(m^2+m+1)\cdot r^{(n)}(\alpha).$$
After simplification we get
$$F'_{n+1}(X)=\frac{n+1}{(X+\alpha+1)^{n+2}}\cdot(\alpha^2+\alpha+1)\cdot\left(f^{(n)}_m(X)\cdot f^{(n)}_m(\alpha)-(m^2+m+1)\cdot r^{(n)}(X)\cdot r^{(n)}(\alpha)\right).$$
The induction hypothesis (\ref{base}) implies
$$F'_{n+1}(X)=\frac{n+1}{(X+\alpha+1)^{n+2}}\cdot(\alpha^2+\alpha+1)\cdot(X+\alpha+1)^n\cdot f^{(n)}_m\left(\frac{\alpha X-1}{X+\alpha+1}\right)=$$
$$=\frac{n+1}{(X+\alpha+1)^{2}}\cdot(\alpha^2+\alpha+1)\cdot f^{(n)}_m\left(\frac{\alpha X-1}{X+\alpha+1}\right)=\left(f^{(n+1)}_m\right)'\left(\frac{\alpha X-1}{X+\alpha+1}\right).$$
In the last step we applied Lemma \ref{deri}, and the derivation rules.
This shows, that 
$$F_{n+1}(X)=f^{(n+1)}_m\left(\frac{\alpha X-1}{X+\alpha+1}\right)+c,$$
with an appropriate constant $c$.
We have to show that $c=0$. In order to prove this, multiply both sides of this equation by $(X+\alpha+1)^{n+1}$:
$$f^{(n+1)}_m(\alpha)\cdot f^{(n+1)}_m(X)-(m^2+m+1)\cdot r^{(n+1)}(\alpha)\cdot r^{(n+1)}(X)=$$
$$(X+\alpha+1)^{n+1}\cdot \left(f^{(n+1)}_m\left(\frac{\alpha X-1}{X+\alpha+1}\right)+c\right)$$
Thus, both sides are polynomials in $X$ of degree $n+1$. The leading coefficient of $X$ in the left hand side is equal to 
$$f^{(n+1)}_m(\alpha)\cdot g(0)-(m^2+m+1)\cdot r^{(n+1)}(\alpha)\cdot h(0)=f^{(n+1)}_m(\alpha).$$
The right hand side is 
$$(X+\alpha+1)^{n+1}\cdot \left(f^{(n+1)}_m\left(\frac{\alpha X-1}{X+\alpha+1}\right)+c\right)=$$
$$\left(\sum_{i=0}^{n+1}\binom{n+1}{i}\cdot(\alpha X-1)^i\cdot(X+\alpha+1)^{n+1-i}\cdot g(n+1-i)\right)+c(X+\alpha+1)^{n+1},$$
and the leading coefficient of $X$ in this polynomial is equal to
$$\left(\sum_{i=0}^{n+1}\binom{n+1}{i}\cdot\alpha^i\cdot g(n+1-i)\right)+c=f^{(n+1)}_m(\alpha)+c.$$
The leading coefficients have to be the same, which implies that $c=0$, and therefore by mathematical induction, the equation (\ref{vtrans}) is true for any $n\in\N$.
\end{prf}

We will use the following interesting lemma about the periodic alternating multisection of sums of binomial coefficients. 

\begin{lemma}\label{remus}
Let 
$$p(i):=\left\lbrace\begin{array}{lr}
a, & \mbox{ if } i\equiv 0 \mod{6},\\
b, & \mbox{ if } i\equiv 1 \mod{6},\\
c, & \mbox{ if } i\equiv 2 \mod{6},\\
-a, & \mbox{ if } i\equiv 3 \mod{6},\\
-b, & \mbox{ if } i\equiv 4 \mod{6},\\
-c, & \mbox{ if } i\equiv 5 \mod{6}.
\end{array}\right.$$
with arbitrary complex numbers $a,b,c$, and let
$$t_n:=\sum_{i=0}^{n}\binom{n}{i}\cdot p(i).$$
Then $$t_{n+12}=3^6\cdot t_n.$$
\end{lemma}

\begin{prf}
We will use the formula of multisections of sums of binomial coefficients (see e.g. \cite{ramus},\cite{ramus2}). Let $0\leq u< v$ be integers, then
$$S_n(u,v):=\binom{n}{u}+\binom{n}{u+v}+\binom{n}{u+2v}+\ldots=\frac{1}{v}\sum_{i=0}^{v-1}\left(2\cos\frac{\pi i}{v}\right)^n\cdot\cos\frac{\pi (n-2u)i}{v}.$$
By this formula,
\begin{equation}\label{tn}
t_n=a\cdot(S_n(0,6)-S_n(3,6))+
b\cdot(S_n(1,6)-S_n(4,6))+
c\cdot(S_n(2,6)-S_n(5,6)).
\end{equation}
The coefficient of $a$ is:
$$S_n(0,6)-S_n(3,6)=\frac{1}{6}\sum_{i=0}^{5}\left(2\cos\frac{\pi i}{6}\right)^n\cdot\left(\cos\frac{\pi n i}{6}-\cos\frac{\pi (n-6)i}{6}\right)=$$
$$=\frac{1}{6}\sum_{i=0}^{5}\left(2\cos\frac{\pi i}{6}\right)^n\cdot\cos\frac{\pi n i}{6}\left(1-\cos\pi i\right).$$
For $i=3$, $\cos\frac{\pi i}{6}=0$;\\
for $i=0,2,4$, $\left(1-\cos\pi i\right)=0$;\\
for $i=1,5$, $\left(1-\cos\pi i\right)=2$, so
$$S_n(0,6)-S_n(3,6)=\frac{1}{3}\left( \left(2\cos\frac{\pi}{6}\right)^n\cdot\cos\frac{\pi n}{6}+\left(2\cos\frac{5\pi}{6}\right)^n\cdot\cos\frac{5\pi n}{6} \right)=$$
$$=\frac{1}{3}\sqrt{3}^n\left(\cos\frac{\pi n}{6}+(-1)^n\cdot\cos\frac{5\pi n}{6} \right).$$
Hence,
$$S_{n+12}(0,6)-S_{n+12}(3,6)=\frac{1}{3}\sqrt{3}^{n+12}\left(\cos\frac{\pi (n+12)}{6}+(-1)^{n+12}\cdot \cos\frac{5\pi (n+12)}{6} \right)=$$
$$=3^6\cdot (S_n(0,6)-S_n(3,6)).$$
Similarly:
$$S_{n+12}(1,6)-S_{n+12}(4,6)=3^6\cdot (S_n(1,6)-S_n(4,6)),$$
$$S_{n+12}(2,6)-S_{n+12}(5,6)=3^6\cdot (S_n(2,6)-S_n(5,6)).$$
Finally by (\ref{tn})
$$t_{n+12}=3^6\cdot t_n.$$
\end{prf}
By this lemma, we are able to calculate $r^{(n)}(\varepsilon_3)$, where $\varepsilon_3$ is a primitive cube root of unity.\\
Since $$r^{(n)}(\varepsilon_3)=\sum_{i=0}^{n}\binom{n}{i}\cdot\varepsilon_3^i\cdot h(n-i),$$
we can use Lemma \ref{remus} with 
$$p(i)=\varepsilon_3^{i}\cdot h(n-i).$$
It is easy to see, that $p(i)$ is of the form we prescribed in Lemma \ref{remus}, since $$p(i+3)=\varepsilon_3^{i+3}\cdot h(n-i-3)=-\varepsilon_3^{i}\cdot h(n-i)=-p(i).$$
Therefore we get
$$r^{(n+12)}(\varepsilon_3)=3^6\cdot r^{(n)}(\varepsilon_3).$$
By calculating $r^{(n)}(\varepsilon_3)$ for $n=1,\ldots,12$, we get
\begin{equation}\label{rne3}
r^{(n)}(\varepsilon_3)=-(-I\sqrt{3})^{n-1}.
\end{equation}
As $3^6=(-I\sqrt{3})^{12}$, (\ref{rne3}) is true for any $n\in\N$.\\
An important consequence is that $\varepsilon_3$ is not a root of $r^{(n)}(X)$.\\
In the following we show, that there is an interesting connection between the roots of $f_m^{(n)}(X)$ and $r^{(n)}(X)$.

\begin{corollary}\label{rootgen}
Let $n\geq2$. If $\beta$ is a root of $f^{(n)}_m(X)$ and $\alpha$ is a root of $r^{(n)}(X)$, then 
$$\frac{\alpha\beta-1}{\beta+\alpha+1}$$
is also a root of $f^{(n)}_m(X)$.
\end{corollary}
\begin{prf}
First we show that if $\beta$ is a root of $f^{(n)}_m(X)$ and  $\alpha$ is a root of $r^{(n)}(X)$, then $\beta+\alpha+1\neq 0$.
We will use, that 
\begin{equation}\label{rnminus}
r^{(n)}(X)=(-1)^{n-1}\cdot r^{(n)}(-X-1).
\end{equation} 
This is trivially true for $n=1$. Assume, that it is true for $n$, then 
$$\left(r^{(n+1)}(X)\right)'=(n+1)\cdot r^{(n)}(X)=(n+1)\cdot(-1)^{n-1}\cdot r^{(n)}(-X-1)=\left((-1)^n\cdot r^{(n+1)}(-X-1)\right)'.$$
The leading coefficients of $r^{(n+1)}(X)$ and $(-1)^n\cdot r^{(n+1)}(-X-1)$ are the same, so (\ref{rnminus}) is true for any $n\geq 1$.\\
The second step is to show that $f^{(n)}_m(X)$ and $r^{(n)}(X)$ can not have a common root.
By Lemma \ref{fr}, if $\alpha$ was a common root of $f^{(n)}_m(X)$ and $r^{(n)}(X)$, and the determinant of 
$$\left(\begin{array}{rr}
\alpha-m & m^2+m+1\\
-1 & \alpha+m+1
\end{array}\right)$$
is not zero, then $\alpha$ would also be a common root of $f^{(n-1)}_m(X)$ and $r^{(n-1)}(X)$. However, the determinant of this matrix is $\alpha^2+\alpha+1$, which is not zero, since  $\alpha$ is a root of $r^{(n)}(X)$, and by (\ref{rne3}), the cube roots of unity are not roots of $r^{(n)}(X)$.\\
Iteratively using this correspondence, we get, that $\alpha$ is a common root of $f^{(1)}_m(X)$ and $r^{(1)}(X)$, but $r^{(1)}(X)=-1$ does not have any root, which is contradiction.\\
If $n\geq 2$ and $\alpha$ is a root of $r^{(n)}(X)$, then $-\alpha-1$ is also a root of $r^{(n)}(X)$. Therefore, as $f^{(n)}_m(X)$ and $r^{(n)}(X)$ don't have a common root, $\beta\neq -\alpha-1$, hence $\beta+\alpha+1\neq0$.\\
Finally, by Theorem \ref{ff}:
$$(\beta+\alpha+1)^n\cdot f^{(n)}_m\left(\frac{\alpha\beta-1}{\beta+\alpha+1}\right)=f^{(n)}_m(\alpha)\cdot f^{(n)}_m(\beta)-(m^2+m+1)\cdot r^{(n)}(\alpha)\cdot r^{(n)}(\beta)=$$
$$=f^{(n)}_m(\alpha)\cdot 0-(m^2+m+1)\cdot 0\cdot r^{(n)}(\beta)=0,$$
whence 
$$\frac{\alpha\beta-1}{\beta+\alpha+1}$$
is also a root of $f^{(n)}_m(X)$.
\end{prf}
Now we will prove, that if $n\geq 2$, then  there exists a root $\alpha$ of $r^{(n)}(X)$, such that the order of the matrix 
$$M=\left(
\begin{array}{rr}
\alpha & -1\\
1 & \alpha+1
\end{array}\right)\in PGL_2(\Q(\alpha))$$
is equal to $n$. Let $\varepsilon_6$ be a primitive sixth root of unity, and $\varepsilon_n$ be a primitive $n$-th root of unity. 
Set
\begin{equation}\label{alfa}
\alpha=\varepsilon_6\cdot\frac{\varepsilon_6+\varepsilon_n}{1-\varepsilon_n}.
\end{equation}
First we show, that $\alpha$ is a root of $r^{(n)}(X)$. 
It is easy to see, that
\begin{equation}\label{aeae}
\frac{\alpha+\varepsilon_6^5}{\alpha+\varepsilon_6}=\varepsilon_n.
\end{equation}
Hence
$$\sum_{i=0}^{n-1}(\alpha+\varepsilon_6^5)^{n-1-i}(\alpha+\varepsilon_6)^i=0,$$
and $\alpha+\varepsilon_6$ is a root of 
$$R^{(n)}(X):= \sum_{i=0}^{n-1}(X+\varepsilon_6^5-\varepsilon_6)^{n-1-i}X^i=\sum_{i=0}^{n-1}(X-I\sqrt{3})^{n-1-i}X^i.$$
We will show that
$$R^{(n)}(X)=-r^{(n)}(X-\varepsilon_6).$$
For $k=0,\ldots,n-1$, the coefficient of $X^k$ in $R^{(n)}(X)$ is:
$$\sum_{i=0}^{k}\binom{n-1-i}{k-i}\cdot(-I\sqrt{3})^{n-1-i-(k-i)}=(-I\sqrt{3})^{n-1-k}\cdot\sum_{i=0}^{k}\binom{n-1-i}{k-i}=$$
$$=(-I\sqrt{3})^{n-1-k}\cdot\binom{n}{k}.$$
The coefficient of $X^k$ in $r^{(n)}(X-\varepsilon_6)$ is:
$$\sum_{i=k}^{n}\binom{n}{i}\cdot\binom{i}{k}\cdot(-\varepsilon_6)^{i-k}\cdot h(n-i)=\sum_{i=k}^{n}\binom{n}{k}\cdot \binom{n-k}{i-k}\cdot(-\varepsilon_6)^{i-k}\cdot h(n-i)=$$
$$=\binom{n}{k}\cdot \sum_{j=0}^{n-k}\binom{n-k}{j}\cdot (-\varepsilon_6)^{j}\cdot h(n-k-j)=\binom{n}{k}\cdot r_{n-k}(-\varepsilon_6)=-\binom{n}{k}\cdot (-I\sqrt{3})^{n-k-1}.$$
In the last step we used equation (\ref{rne3}), and that $-\varepsilon_6$ is a primitive cube root of unity.
Since all coefficients of $R^{(n)}(X)$ are the negative of the corresponding coefficients of $r^{(n)}(X-\varepsilon_6)$, we obtain $$R^{(n)}(X)=-r^{(n)}(X-\varepsilon_6).$$
Thus $\alpha+\varepsilon_6$ is a root of $r^{(n)}(X-\varepsilon_6)$, and therefore $\alpha$ is a root of $r^{(n)}(X)$.\\
Due to the special form of $\alpha$, it is also easy to determine the order of 
$$M=\left(
\begin{array}{rr}
\alpha & -1\\
1 & \alpha+1
\end{array}\right)\in PGL_2(\Q(\alpha)).$$
The eigenvalues of this matrix are $\alpha+\varepsilon_6$ and $\alpha+\varepsilon_6^5$, and by \eqref{aeae} the quotient of the eigenvalues of $M$ is $\varepsilon_n$, therefore $M\in PGL_2(\Q(\alpha))$ is a matrix of order $n$.\\
This allows us to determine the Galois group of $f^{(n)}_m(X)$ over $\Q(\alpha)$.
\begin{theorem}
Let $n\geq2$. If $\alpha$ is defined as in (\ref{alfa}), and $m$ is a rational number, such that $f^{(n)}_m(X)$ is irreducible over $\Q(\alpha)$, then the Galois group of $f^{(n)}_m(X)$ over $\Q(\alpha)$ is the cyclic group of order $n$.
\end{theorem}
\begin{prf}
By Corollary \ref{rootgen}, if $\beta$ is a root of $f^{(n)}_m(X)$, then 
$$\frac{\alpha\beta-1}{\beta+\alpha+1}$$
is also a root of $f^{(n)}_m(X)$. As the order of the matrix 
$$M=\left(
\begin{array}{rr}
\alpha & -1\\
1 & \alpha+1
\end{array}\right)\in PGL_2(\Q(\alpha))$$
is $n$,  $\{ \sigma^i(\beta) | i=0,\ldots,n-1 \}$ is the complete set of roots of $f^{(n)}_m(X)$, where
$$\sigma:z\mapsto\frac{\alpha z-1}{z+\alpha+1}.$$
Finally $\sigma^i(\beta)\in \Q(\alpha,\beta)$ for all $i=0,\ldots,n-1$, hence the splitting field of $f^{(n)}_m(X)$ over $\Q(\alpha)$ is $\Q(\alpha,\beta)$. Therefore, $\Q(\alpha,\beta)/\Q(\alpha)$, is a cyclic Galois extension of degree $n$.
\end{prf}
\subsection{The discriminant of $f^{(n)}_m(X)$}\label{Sdisc}

In this section, we determine the discriminants of the polynomials $f^{(n)}_m(X)$.\\
We will use the following properties of the resultant $\res(A,B)$ of polynomials 
$$A(X)=a_{n}X^n+a_{n-1}X^{n-1}+\ldots+a_1X+a_0,$$
$$B(X)=b_{m}X^m+b_{m-1}X^{m-1}+\ldots+b_1X+b_0.$$
(see for e.g. Chapter 2.3.3. \cite{pohst}).
\begin{equation}\label{csere}
\res(A,B)=(-1)^{nm}\cdot\res(B,A).
\end{equation}
If $n\geq m$ and $Q,R$ are polynomials such that $A=QB+R$, where $\deg(R)=r$, then
\begin{equation}\label{longdiv}
\res(A,B)=b_m^{n-r}\cdot\res(R,B).
\end{equation}
For any monic polynomial $Q$,
\begin{equation}\label{szor}
\res(AQ,B)=\res(A,B)\cdot\res(Q,B).
\end{equation}
If $\lambda\in\R$, then
\begin{eqnarray}\label{scalarmult}
\res(\lambda\cdot A,B)=\lambda^m\cdot\res(A,B),\\
\res(A,\lambda\cdot B)=\lambda^n\cdot\res(A,B).
\end{eqnarray}

\begin{lemma}\label{resff}
For $n\geq2$ the resultant of $f^{(n)}_m(X)$ and $f^{(n-1)}_m(X)$ with respect to $X$ is
$$\res\left(f^{(n)}_m(X),f^{(n-1)}_m(X)\right)=(m^2+m+1)^{n-1}\cdot3^{\frac{(n-1)(n-2)}{2}}\cdot(-1)^{\frac{n(n-1)}{2}}.$$
\end{lemma}

\begin{prf}
First we determine the resultant of $r^{(n)}(X)$ and $X^2+X+1$.\\
The remainder of $r^{(n)}(X)$ divided by $X^2+X+1$ is 
$$T_n(X)=\sum_{i=0}^{n}\binom{n}{i}\cdot G(i)\cdot h(n-i),$$
where 
$$G(i)=\left\lbrace
\begin{array}{lr}
1, & \mbox{ if }i\equiv0\mod{3},\\
X, & \mbox{ if }i\equiv1\mod{3},\\
-X-1, & \mbox{ if }i\equiv2\mod{3}.
\end{array}\right.$$
Thus, we can write $T_n(X)=A+B\cdot X$, with
$$A=\sum_{i=0}^{n}\binom{n}{i}\cdot G_0(i)\cdot h(n-i),$$
$$B=\sum_{i=0}^{n}\binom{n}{i}\cdot G_1(i)\cdot h(n-i),$$
$$G_0(i):=\left\lbrace\begin{array}{lr}
1, & \mbox{ if } i\equiv 0 \mod{3},\\
0, & \mbox{ if } i\equiv 1 \mod{3},\\
-1, & \mbox{ if } i\equiv 2 \mod{3}.\\
\end{array}\right.\qquad\mbox{and}\qquad
G_1(i):=\left\lbrace\begin{array}{lr}
0, & \mbox{ if } i\equiv 0 \mod{3},\\
1, & \mbox{ if } i\equiv 1 \mod{3},\\
-1, & \mbox{ if } i\equiv 2 \mod{3}.\\
\end{array}\right.$$
By applying Lemma \ref{remus} to the coefficients of $T_n(X)$, we obtain, $T_{n+12}(X)=3^6\cdot T_n(X)$, whence
using the identities (\ref{longdiv}) and  (\ref{scalarmult}),
$$\res\left(r^{(n+12)}(X),X^2+X+1\right)=\res\left(T_{n+12}(X),X^2+X+1\right)=$$
$$\res\left(3^6\cdot T_n(X),X^2+X+1\right)=3^{12}\cdot \res\left(T_n(X),X^2+X+1\right)=$$
$$3^{12}\cdot \res\left(r^{(n)}(X),X^2+X+1\right).$$
By calculating $\res\left(r^{(n)}(X),X^2+X+1\right)$ for $n=1,\ldots,12$, we get that 
$$\res\left(r^{(n)}(X),X^2+X+1\right)=3^{n-1},$$
for any $n\geq 1$, and by (\ref{csere}), 
\begin{equation}\label{resrnxx}
\res\left(X^2+X+1,r^{(n)}(X)\right)=(-1)^{2(n-1)}\cdot3^{n-1}=3^{n-1}.
\end{equation}
The next step is to determine the resultant of $r^{(n)}(X)$ and $f^{(n)}_m(X)$ with respect to $X$.\\
By (\ref{fn+1}) and (\ref{rn+1}), we have
$$f^{(n)}_m(X)=-(X-m)\cdot r^{(n)}(X)+(X^2+X+1)\cdot r^{(n-1)}(X).$$
Therefore, by (\ref{csere}), (\ref{longdiv}), (\ref{szor}) and (\ref{resrnxx}),
$$\res\left(r^{(n)}(X),f^{(n)}_m(X)\right)=(-1)^{n(n-1)}\cdot\res\left(f^{(n)}_m(X),r^{(n)}(X)\right)=$$
$$n^0\cdot\res\left((X^2+X+1)\cdot r^{(n-1)}(X),r^{(n)}(X)\right)=$$
$$\res\left(X^2+X+1,r^{(n)}(X)\right)\cdot\res\left(r^{(n-1)}(X),r^{(n)}(X)\right)=3^{n-1}\res\left(r^{(n-1)}(X),r^{(n)}(X)\right).$$
By (\ref{rn+1}), we have
$$r^{(n)}(X)=(X+m+1)\cdot r^{(n-1)}(X)-f^{(n-1)}_m(X),$$ whence by (\ref{longdiv}), (\ref{szor}) and (\ref{csere})
$$\res\left(r^{(n)}(X),r^{(n-1)}(X)\right)=(n-1)^0\cdot\res\left(-f^{(n-1)}_m(X),r^{(n-1)}(X)\right)=$$
$$(-1)^{n-2}\cdot\res\left(f^{(n-1)}_m(X),r^{(n-1)}(X)\right)=(-1)^{(n-1)(n-2)+n-2}\cdot\res\left(r^{(n-1)}(X),f^{(n-1)}_m(X)\right)=$$
$$(-1)^n\cdot\res\left(r^{(n-1)}(X),f^{(n-1)}_m(X)\right).$$
Hence, 
$$\res\left(r^{(n)}(X),f^{(n)}_m(X)\right)=(-1)^{n}\cdot 3^{n-1}\cdot\res\left(r^{(n-1)}(X),f^{(n-1)}_m(X)\right).$$
By iterating this equation, and using that $\res(r^{(1)}(X),f^{(1)}_m(X))=-1$, we obtain
\begin{equation}\label{resfr}
\res\left(r^{(n)}(X),f^{(n)}_m(X)\right)=3^{\frac{n(n-1)}{2}}\cdot(-1)^{\frac{n(n+1)}{2}}.
\end{equation}
Finally, we obtain the resultant of $f^{(n)}_m(X)$ and $f^{(n-1)}_m(X)$ with respect to $X$.\\
By (\ref{fn+1}), we have
$$f^{(n)}_m(X)=(X-m)\cdot f^{(n-1)}_m(X)+(m^2+m+1)\cdot r^{(n-1)}(X).$$
Therefore, by (\ref{longdiv}), (\ref{scalarmult}) and (\ref{resfr})
$$\res\left(f^{(n)}_m(X),f^{(n-1)}_m(X)\right)=\res\left((m^2+m+1)\cdot r^{(n-1)}(X),f^{(n-1)}_m(X)\right)=$$
$$(m^2+m+1)^{n-1}\cdot \res\left(r^{(n-1)}(X),f^{(n-1)}_m(X)\right)=$$
$$(m^2+m+1)^{n-1}\cdot3^{\frac{(n-1)(n-2)}{2}}\cdot(-1)^{\frac{n(n-1)}{2}}.$$
\end{prf}

\begin{theorem}
For $n\geq 2$ the discriminant of $f^{(n)}_m(X)$ with respect to $X$ is 
\begin{equation}\label{discrim}
D_{f^{(n)}_m}=3^{\frac{(n-1)(n-2)}{2}}\cdot n^n\cdot (m^2+m+1)^{n-1}.
\end{equation}
\end{theorem}
\begin{prf}
The discriminant of $f^{(n)}_m(X)$ with respect to $X$, is 
$$D_{f^{(n)}_m}=(-1)^{\frac{n(n-1)}{2}}\cdot \res\left(f^{(n)}_m(X),\left(f^{(n)}_m\right)'(X)\right).$$
By Lemma \ref{deri} and (\ref{scalarmult}):
$$D_{f^{(n)}_m}=(-1)^{\frac{n(n-1)}{2}}\cdot \res\left(f^{(n)}_m(X),n\cdot f^{(n-1)}_m(X)\right)=n^n\cdot(-1)^{\frac{n(n-1)}{2}}\cdot \res\left(f^{(n)}_m(X),f^{(n-1)}_m(X)\right).$$
Thus by Lemma \ref{resff}, 
$$D_{f^{(n)}_m}=n^n\cdot(m^2+m+1)^{n-1}\cdot(-1)^{\frac{n(n-1)}{2}}\cdot \res\left(r^{(n-1)}(X),f^{(n-1)}_m(X)\right)=$$
$$=3^{\frac{(n-1)(n-2)}{2}}\cdot n^n\cdot(m^2+m+1)^{n-1}.$$
\end{prf}
\subsection{Irreducibility of $f^{(n)}_m(X)$}\label{Sired}
In this section, we will give a sufficient condition for the irreducibility of $f^{(n)}_m(X)$ over $\Q$. For the investigation of an integral basis of these polynomials, we will use exactly the same condition, and therefore we do not need the necessary condition of the irreducibility.\\
Here $v_p(x)$ denotes the $p$-adic valuation of a rational number $x$.

\begin{lemma}\label{eisen}
Let $n\geq 2$. If there exist a prime $p\neq3$, such that $v_p(m^2+m+1)=1$, then $v_p\left(f^{(n)}_m(m)\right)=1$.
\end{lemma}
\begin{prf}
By Lemma \ref{fr}, if $n\geq2$ then 
$$r^{(n)}(m)=(2m+1)\cdot r^{(n-1)}(m)-f^{(n-1)}_m(m)$$
and
$$f^{(n)}_m(m)=(m^2+m+1)\cdot r^{(n-1)}(m),$$
thus
$$
r^{(n)}(m)=(2m+1)\cdot r^{(n-1)}(m)-(m^2+m+1)\cdot r^{(n-2)}(m).$$
We have $v_p(2m+1)=0$, since $4(m^2+m+1)-(2m+1)^2=3$ and $p\neq 3$.\\
Now assume that $v_p\left(r^{(n-1)}(m)\right)=v_p\left(r^{(n-2)}(m)\right)$, then
$$v_p\left((2m+1)\cdot r^{(n-1)}(m)\right)=0+v_p\left(r^{(n-1)}(m)\right),$$
$$v_p\left((m^2+m+1)\cdot r^{(n-2)}(m)\right)=1+v_p\left(r^{(n-1)}(m)\right),$$
whence
$$v_p\left(r^{(n)}(m)\right)=\min(v_p\left(r^{(n-1)}(m)\right),v_p\left(r^{(n-1)}(m)\right)+1)=v_p\left(r^{(n-1)}(m)\right).$$
Therefore, $v_p\left(r^{(1)}(m)\right)=v_p(-1)=0$ and $v_p\left(r^{(2)}(m)\right)=v_p(-2m-1)=0$ implies
$$v_p\left(r^{(n)}(m)\right)=0$$
and 
$$v_p\left(f^{(n)}_m(m)\right)=v_p\left(m^2+m+1\right)+v_p\left(r^{(n-1)}(m)\right)=1.$$
\end{prf}

\begin{corollary} \label{irreduc}
Let $n\geq 2$. If there exist a prime $p\neq3$, such that $v_p(m^2+m+1)=1$, then $f^{(n)}_m(X)$ is irreducible over $\Q$.
\end{corollary}
\begin{prf}
The Taylor series of $f^{(n)}_m(X)$ at $m$ is:
$$f^{(n)}_m(X)=\sum_{i=0}^{n}\frac{\left({f^{(n)}_m}\right)^{(i)}(m)}{i!}\cdot(X-m)^i.$$
Thus, by shifting the polynomial by $m$:
$$f^{(n)}_m(X+m)=\sum_{i=0}^{n}\frac{\left({f^{(n)}_m}\right)^{(i)}(m)}{i!}\cdot X^i.$$
By Lemma \ref{deri}
$$\left({f^{(n)}_m}\right)^{(i)}(X)=\frac{n!}{(n-i)!}\cdot f^{(n-i)}_m(X),$$
whence
$$f^{(n)}_m(X+m)=\sum_{i=0}^{n}\binom{n}{i}\cdot f^{(n-i)}_m(m)\cdot X^i.$$
With $f^{(0)}_m(m)=1$ and $f^{(1)}_m(m)=0$, Lemma \ref{eisen} shows that if $n\geq 2$, then $f^{(n)}_m(X+m)$ is Eisenstein with respect to $p$, so it is irreducible, and consequently $f^{(n)}_m(X)$ is also irreducible.
\end{prf}

\section{Periodically repeating integral bases}\label{prib}
In this section we will prove, that an integral basis of the generalized simplest fields generated by a root of $f^{(n)}_m(X)$ $(n\geq2)$ is repeating periodically, and we give an upper bound for the period length $n_0$, in terms of $n$.\\

The periodic property of an integral basis makes sense if we investigate the integral bases of a parametric family of number fields of type $\Q(\beta_t)$, where $\beta_t$ is a root of a polynomial $f^{(n)}_{t}(X)\in\Z[t][X]$ of degree $n$, with integer parameter $t\in\Z$. In this case we determine an integral basis of $\Q(\beta_t)$ of the form 
$$\left(h_0(\beta_t),h_1(\beta_t),h_2(\beta_t),\ldots,h_{n-1}(\beta_t)\right), $$
where $h_i(X)\in\Q[X]$, $(i=0,\ldots,n-1)$, are polynomials with $\deg(h_i(X))=i$.

In order to construct a parametric family of polynomials with integer coefficients from $f_m^{(n)}$, let
\begin{itemize}
\item $m=t$, if $n\equiv1,2\mod{3},$
\item $m=\frac{t}{3}$, if $n\equiv0\mod{3}$.
\end{itemize}
It is easy to check, that for any $t\in\Z$, the generalized simplest polynomials $f^{(n)}_m(X)$ have integer coefficients. In the following we write $f^{(n)}_t(X)$, and it refers to the definition above.\\
In this case the modified results of the previous section are the following.\\\\
The discriminant of $f_t^{(n)}(X)$ with respect to $X$ is
\begin{equation}\label{tdisc}
D_{f_t^{(n)}}=
\left\lbrace\begin{array}{lr}
3^{\frac{(n-1)(n-2)}{2}}\cdot n^n\cdot (t^2+t+1)^{n-1},& $if $ n\equiv1,2\mod{3},\\
3^{\frac{(n-1)(n-6)}{2}}\cdot n^n\cdot (t^2+3t+9)^{n-1},& $if $n\equiv0\mod{3}.
\end{array}\right.
\end{equation}
If $t\in\Z$ is an integer parameter, such that there exists a prime $p\neq3$, for which 
\begin{itemize}
\item $v_p(t^2+t+1)=1$, if $n\equiv1,2\mod{3},$
\item $v_p(t^2+3t+9)=1$, if $n\equiv0\mod{3},$
\end{itemize}
then $f_t^{(n)}(X)$ is irreducible over $\Q$.\\\\
Let $\beta_t$ be a root of $f^{(n)}_t(X)$. We recall that an integral basis of the fields $\Q(\beta_t)$ is \textit{repeating periodically modulo $n_0$}, if for each residue class $r$ modulo $n_0$, there exist polynomials $h^{(r)}_i(X)\in\Q[X]$, $(i=0,\ldots n-1)$, such that, if $t\equiv r\mod{n_0}$, and $f^{(n)}_t(X)$ is irreducible then  
$$\left(h^{(r)}_0(\beta_t),h^{(r)}_1(\beta_t),h^{(r)}_2(\beta_t),\ldots,h^{(r)}_{n-1}(\beta_t)\right), $$
is an integral basis of $\Q(\beta_t)$. (See \cite{ibm},\cite{sextic},\cite{ibpf}).\\
There is an equivalent definition of a periodically repeating integral basis which is more convenient to use in practice.
\begin{lemma}\label{ekvparam}
Let $f_t^{(n)}(X)\in\Z[t,X]$ be a parametric family of polynomials, where $t$ is an integer parameter such that $f_t^{(n)}(X)$ is irreducible, and let $\beta_t$ be a root of $f_t^{(n)}(X)$.\\
An integral basis of the fields $\Q(\beta_t)$ is repeating periodically modulo $n_0$, if and only if for any $k(X)\in\Q[X]$, $\deg(k(X))\leq n-1$, the following property is satisfied:\\
for any parameters $t,s\in\Z$ for which $f_{t}^{(n)}(X)$ and $f_{s}^{(n)}(X)$ are irreducible and $n_0\mid(t-s)$, 
$k(\beta_{t})$ is an algebraic integer if and only if $k(\beta_{s})$ is an algebraic integer.
\end{lemma}

\begin{prf}
First assume, that an integral basis of $\Q(\beta_t)$ is repeating periodically modulo $n_0$. Let $V=\sum_{i=0}^{n-1}\Q X^i$. Then the $\Q$-linear map $\varphi:X^i\mapsto\beta_t^i$ gives an isomorphism of the $\Q$-vector spaces $V$ and $\Q(\beta_t)$ for any $t\in\N$. Let $M$ be the $\Z$-module generated by 
$$h^{(r)}_0(X),h^{(r)}_1(X),h^{(r)}_2(X),\ldots,h^{(r)}_{n-1}(X),$$
where $t\equiv r\mod{n_0}$. If we restrict $\varphi$ to $M$ then we obtain an isomorphism of $M$ and the ring of integers of $\Q(\beta_t)$ for any $t\equiv r\mod{n_0}$. This means that if $t\equiv r\mod{n_0}$, then for any $k(X)\in\Q[X]$ with $\deg(k(X))\leq n-1$, $k(\beta_t)$ is an algebraic integer if and only if $k(X)\in M$. Therefore if $t\equiv s\equiv r\mod{n_0}$, then $k(\beta_{t})$ is an algebraic integer if and only if $k(\beta_{s})$ is an algebraic integer.\\
On the other hand, assume that for any polynomial $k(X)\in\Q[X]$ with $\deg(k(X))\leq n-1$, $k(\beta_t)$ is an algebraic integer if and only if $k(\beta_s)$ is an algebraic integer. Let $r\equiv t\equiv s\mod{n_0}$ and $h^{(r)}_i(X)\in\Q[X]$, $(i=0,\ldots n-1)$ be polynomials, such that 
$$\left(h^{(r)}_0(\beta_t),h^{(r)}_1(\beta_t),h^{(r)}_2(\beta_t),\ldots,h^{(r)}_{n-1}(\beta_t)\right)$$
is an integral basis of $\Q(\beta_t)$, where $\deg(h^{(r)}_i(X))=i$ (one can use the Hermite normal form representation of any integral basis of $\Q(\beta_t)$ with respect to the $\Q$-basis $\left\lbrace 1,\beta_t,\beta_t^2,\ldots,\beta_t^{n-1}\right\rbrace$ to obtain such $h^{(r)}_i(X)$ polynomials). Since all elements of an integral basis are algebraic integers, by our assumption we have that $h^{(r)}_i(\beta_s)$, is algebraic integer for any $i=0,\ldots n-1$. We should prove, that the elements  $h^{(r)}_i(\beta_s)$ $(i=0,\ldots n-1)$ form an integral basis of $\Q(\beta_s)$.\\
Let $k(X)\in\Q[X]$, $\deg(k(X))\leq n-1$, be a polynomial, such that $k(\beta_s)$ is an algebraic integer. By our assumption $k(\beta_t)$ is also an algebraic integer, thus we can represent it in the integral basis 
$$\left(h^{(r)}_0(\beta_t),h^{(r)}_1(\beta_t),h^{(r)}_2(\beta_t),\ldots,h^{(r)}_{n-1}(\beta_t)\right)$$
with integer coefficients. Therefore we can write
$$k(X)=k_0\cdot h^{(r)}_0(X)+k_1\cdot h^{(r)}_1(X)+\ldots+k_{n-1}\cdot h^{(r)}_{n-1}(X),\qquad k_0,k_1,\ldots,k_{n-1}\in\Z.$$
By substituting $\beta_s$ in place of $\beta_t$ in the equation above, we obtain that we can represent $k(\beta_s)$ as an integer linear combination of the elements $h^{(r)}_i(\beta_s)$, $(i=0,\ldots n-1)$. Therefore, if these elements are linearly independent over $\Q$, then 
$$\left(h^{(r)}_0(\beta_s),h^{(r)}_1(\beta_s),h^{(r)}_2(\beta_s),\ldots,h^{(r)}_{n-1}(\beta_s)\right)$$
is an integral basis of $\Q(\beta_s)$, and an integral basis of the fields $\Q(\beta_t)$ is repeating periodically modulo $n_0$. However, it is trivial, since $\deg(h^{(r)}_i(X))=i$, $(i=0,\ldots n-1)$, and the minimal polynomial of $\beta_s$ is the irreducible polynomial $f_s^{(n)}(X)$ of degree $n$.
\end{prf}

\subsection{Upper bound for the period length}

In this section we prove, that if the 3-free part of $t^2+t+1$ (or $t^2+3t+9$ respectively) is square-free, then an integral basis of $\Q(\beta_t)$ is repeating periodically modulo $n_0$, where
$$n_0^2\mid\left(3^{\frac{n^2}{2}}\cdot n^n\right)^n.$$
First we show, that if the 3-free part of $t^2+t+1$ (or $t^2+3t+9$ respectively) is square-free, then the index $\ind(\beta_t)$ of $\beta_t$ is independent of $t$. Therefore, for any $n\in\N$ there exists an integer $C_n$, such that for any algebraic integer $\gamma_t\in\Q(\beta_t)$:
$$C_n\cdot\gamma_t\in \Z[\beta_t].$$
It means, that we can represent $\gamma_t$ in the $\Q$-basis 
$$\left\lbrace 1,\beta_t,\beta_t^2,\ldots,\beta_t^{n-1}\right\rbrace,$$
with rational coefficients having denominator $d\mid C_n$.\\
Finally, we prove that if such $C_n$ exists, then an integral basis of these fields is repeating periodically modulo $n_0=C_n^n$.\\

The \textit{index} of an algebraic integer $\beta_t$ is defined by
$$\ind(\beta_t):=({\cal O}_{\Q(\beta_t)}^+:\Z[\beta_t]^+),$$
where ${\cal O}_{\Q(\beta_t)}$ is the ring of integers of $\Q(\beta_t)$. That is the index of the additive group of the polynomial ring $\Z[\beta_t]$ in the additive group of ${\cal O}_{\Q(\beta_t)}$. (for more details see for e.g. I.Gaál \cite{book}, Section 1.2).\\
This is in fact, the product of the denominators of the elements of an integral basis of $\Q(\beta_t)$, represented in the $\Q$-basis $\left(1,\beta_t,\ldots,\beta_t^{n-1}\right)$.\\
If $D(\beta_t)$ denotes the discriminant of the minimal polynomial of $\beta_t$ (which is equal to the discriminant of the basis $\left(1,\beta_t,\ldots,\beta_t^{n-1}\right)$), and $D_M$ denotes the discriminant of the field $M=\Q(\beta_t)$, then
\begin{equation}\label{ind2}
D(\beta_t)=\ind(\beta_t)^2\cdot D_M.
\end{equation}
\begin{lemma}\label{ein}
(\cite{Nark} Chapter II., Lemma 2.17)
Let $M$ be a number field of degree $n$, and let $\beta\in M$ be a nonzero algebraic integer of degree $n$. Suppose that the minimal polynomial of $\beta$ is Eisenstein with respect to the prime $p$. Then $p\nmid\ind(\beta_t)$.
\end{lemma}
By this lemma, we can exclude certain prime divisors of the discriminant $D(\beta_t)$.
\begin{theorem}\label{pind}
Let $t$ be an integer and $\beta_t$ be a root of $f^{(n)}_t(X)$. If $p\neq 3$ is a prime, for which 
$$v_p(t^2+t+1)=1,\qquad \mbox{ if } n\equiv1,2\mod{3},$$
$$v_p(t^2+3t+9)=1,\qquad \mbox{ if } n\equiv0\mod{3},$$ 
then $p\nmid\ind(\beta_t)$.
\end{theorem}
\begin{prf}
If $n\equiv1$ or $2\mod{3}$, then by Corollary \ref{irreduc}, the minimal polynomial of $\beta_t-t$ is Eisenstein with respect to $p$, so by Lemma \ref{ein}, $p\nmid\ind(\beta_t-t)$. Since the index of an algebraic integer is invariant under translation, we obtain $p\nmid\ind(\beta_t)$.\\
If $n\equiv0\mod{3}$, then by Corollary \ref{irreduc}, the minimal polynomial of $3\beta_t-t$ is Eisenstein with respect to $p$, so by Lemma \ref{ein}, $p\nmid\ind(3\beta_t-t)$, and therefore $p\nmid\ind(3\beta_t)$. If $p\neq3$, then $v_p(\ind(3\beta_t))=v_p(\ind(\beta_t))$, and thus we have $p\nmid\ind(\beta_t)$ in this case, too.
\end{prf}
By (\ref{tdisc}) and (\ref{ind2})
$$\ind(\beta_t)^2\mid 3^{\frac{(n-1)(n-2)}{2}}\cdot n^n\cdot (t^2+t+1)^{n-1},\qquad \mbox{ if } n\equiv1,2\mod{3},$$
$$\ind(\beta_t)^2\mid 3^{\frac{(n-1)(n-6)}{2}}\cdot n^n\cdot (t^2+3t+9)^{n-1},\qquad \mbox{ if } n\equiv0\mod{3}.$$
Hence, by Theorem \ref{pind}, if the 3-free part of $t^2+t+1$ (or $t^2+3t+9$, respectively) is square-free, then 
$$\ind(\beta_t)^2\mid 3^{\frac{(n-1)(n-2)}{2}+v_3(t^2+t+1)}\cdot n^n,\qquad \mbox{ if } n\equiv1,2\mod{3},$$
$$\ind(\beta_t)^2\mid 3^{\frac{(n-1)(n-6)}{2}+v_3(t^2+3t+9)}\cdot n^n,\qquad \mbox{ if } n\equiv0\mod{3}.$$
But it is easy to check, that for any $t\in\Z$, $v_3(t^2+t+1)\leq1$ and  $v_3(t^2+3t+9)\leq3$.
Thus we get the following corollary.
\begin{corollary}\label{denom}
Let $\beta_t$ be a root of $f^{(n)}_t(X)$, and $C_n$ be the greatest integer, such that
$$C_n^2\mid 3^{\frac{n^2-3n+4}{2}}\cdot n^n,\qquad \mbox{ if } n\equiv1,2\mod{3},$$
$$C_n^2\mid 3^{\frac{n^2-7n+12}{2}}\cdot n^n,\qquad \mbox{ if } n\equiv0\mod{3}.$$
If $t^2+t+1$ (or $t^2+3t+9$, respectively) is square-free, then for any algebraic integer $\gamma_t\in\Q(\beta_t)$,
$$C_n\cdot \gamma_t\in\Z[\beta_t].$$
\end{corollary}
Therefore, in both cases
$$C_n^2\mid 3^{\frac{n^2}{2}}\cdot n^n.$$
Following similar approach as in Theorem 2 of \cite{ibm}, we show, that it implies that an integral basis of these fields is repeating periodically.\\
\begin{lemma}\label{cnn0}
Let $\beta_t$ be a root of $f_t^{(n)}(X)$. If there exists an integer $C_n$, such that for any algebraic integer $\gamma_t\in\Q(\beta_t)$,
$$C_n\cdot\gamma_t\in \Z[\beta_t],$$
then an integral basis of the fields $\Q(\beta_t)$ is repeating periodically modulo $n_0=C_n^n$.
\end{lemma}
\begin{prf}
Let $\gamma_t\in\Q(\beta_t)$ be an algebraic integer and represent it the $\Q$-basis $\left\lbrace 1,\beta_t,\beta_t^2,\ldots,\beta_t^{n-1}\right\rbrace$:
$$\gamma_t=\frac{a_0+a_1\beta_t+a_2\beta_t^2+\ldots+a_{n-1}\beta_t^{n-1}}{d}, \qquad a_i\in\Z, d\mid C_n.$$
Denote the conjugates of $\gamma_t$ by $\gamma_t^{(j)}$, $j=1,\ldots,n$. The defining polynomial of $\gamma_t$ is 
$$\prod_{j=1}^{n}(X-\gamma_t^{(j)})=\frac{1}{d^n}\prod_{j=1}^{n}(dX-a_0-a_1\beta_t^{(j)}-a_2(\beta_t^{(j)})^2-\ldots-a_{n-1}(\beta_t^{(j)})^{n-1}).$$
The product is a symmetric polynomial of $\beta_t^{(1)},\beta_t^{(2)},\ldots,\beta_t^{(n)}$, hence its coefficients can be expressed as integer polynomials of the coefficients of the defining polynomial of $\beta_t$, that is $f_t^{(n)}(X)$. The coefficients of $f_t^{(n)}(X)$ are integer polynomials of $t$, hence there exist polynomials $P_0,\ldots,P_{n-1}\in\Z[X]$, such that
$$\prod_{j=1}^{n}(X-\gamma_t^{(j)})=\frac{1}{d^n}\left(d^nX^n+P_{n-1}(t)(dX)^{n-1}+\ldots+P_1(t)(dX)+P_0(t)\right).$$
Therefore the element $\gamma_t$ is an algebraic integer if $d^n\mid d^jP_j(t)$ for all $j=1,\ldots,n-1$.\\
Now let $s\in\Z$ be an integer, such that  $C_n^n\mid (t-s)$. Then $C_n^n\mid P_j(t)-P_j(s)$, and therefore
$d^n\mid d^jP_j(t)$ if and only if $d^n\mid d^jP_j(s)$.
It implies that for any $k(X)\in\Q[X]$, $(\deg(k(X))\leq n-1)$, $\gamma_t=k(\beta_t)$ is algebraic integer if and only if $\gamma_s=k(\beta_{s})$ is an algebraic integer, thus by Lemma \ref{ekvparam} an integral basis of $\Q(\beta_t)$ is repeating periodically modulo $n_0=C_n^n$.
\end{prf}
By Corollary \ref{denom} and Lemma \ref{cnn0}, we obtain the main result of this section.
\begin{theorem}\label{period}
Let $\beta_t$ be a root of $f^{(n)}_t(X)$, and
$n_0$ be a greatest integer, such that 
$$n_0^2\mid \left(3^{\frac{n^2}{2}}\cdot n^n\right)^n.$$
If $t^2+t+1$ (or $t^2+3t+9$, respectively) is square-free, then, an integral basis of the fields $\Q(\beta_t)$ is repeating periodically modulo $n_0$.
\end{theorem}

\subsection{Special cases}
In this section we give better upper bounds for the period lengths of the integral bases of the fields $\Q(\beta_t)$ of degree $n\leq 12$.\\
We will obtain the better upper bounds for the period length by using the dual bases of certain rational bases of these fields. For the advantage of the dual basis see for e.g. \cite{pohst}.\\\\
Let $M$ be an algebraic number field of degree $n$. If $(\beta_0,\beta_1,\ldots,\beta_{n-1})$ is a $\Q$-basis of $M$, then its dual basis is $(\gamma_0,\gamma_1,\ldots,\gamma_{n-1})$, where $\gamma_i\in M$, such that 
$$Tr_{M/\Q}(\gamma_i\cdot \beta_j)=\delta_{ij}=\left\lbrace
\begin{array}{ll}
1, & $if $i=j,\\
0, & $otherwise.$
\end{array}
\right.
$$
The dual basis of any $\Q$-basis of $M$ always uniquely exists.\\
Let $\left(\beta_0,\beta_1,\ldots,\beta_{n-1}\right)$ be a $\Q$-basis of $M$ containing algebraic integers, and $\left(\gamma_0,\gamma_1,\ldots,\gamma_{n-1}\right)$ its dual basis. Then it is easy to see, that we can represent any algebraic integer $\gamma\in M$ in the $\Q$-basis $\left(\gamma_0,\gamma_1,\ldots,\gamma_{n-1}\right)$ with integer coefficients.\\
Let $\gamma\in M$ be an algebraic integer, and represent is in the $\Q$-basis $\left(\gamma_0,\gamma_1,\ldots,\gamma_{n-1}\right)$:
$$\gamma=f_0\gamma_0,+f_1\gamma_1+\ldots+f_{n-1}\gamma_{n-1}.$$
By the definition of the dual basis,  $f_i=Tr_{M/\Q}(\gamma\cdot \beta_i)$, which is a rational integer, since $\gamma$ and $\beta_i$ both algebraic integers.\\
Let $\beta_t$ be a root of $f^{(n)}_t(X)$, and $\left(\gamma_0,\gamma_1,\ldots,\gamma_{n-1}\right)$ be the dual basis of $\left(1,\beta_t,\beta_t^2\ldots,\beta_t^{n-1}\right)$ in $\Q(\beta_t)$. Represent the elements $\gamma_i$ of the dual basis in the $\Q$-basis $\left(1,\beta_t,\beta_t^2\ldots,\beta_t^{n-1}\right)$:
$$\gamma_i=\sum_{j=0}^{n-1}c_{ij}\beta_t^j.$$
Let $d$ be the least common multiple of the denominators of the rational numbers $c_{ij}$, $(i,j=0,\ldots,n-1)$. Then we can represent any algebraic integer $\gamma\in\Q(\beta_t)$ in the $\Q$-basis $\left(\gamma_0,\gamma_1,\ldots,\gamma_{n-1}\right)$ with integer coefficients, thus in the $\Q$-basis $\left(1,\beta_t,\beta_t^2\ldots,\beta_t^{n-1}\right)$, we can represent $\gamma$ with rational coefficients having denominator $d$. By the results of the previous section, an integral basis of the fields $\Q(\beta_t)$ is repeating periodically modulo $d^n$.\\
We explicitly calculated these dual bases for $n=2,\ldots,12$, and obtained Theorem \ref{periodl}.\\
For instance, if $n=4$, the dual basis of $\left(1,\beta_t,\beta_t^2,\beta_t^3\right)$ is 

$$\displaystyle\left(
\begin{array}{rrrr}
\frac{8t^2+3t+12}{12(t^2+t+1)} & +\frac{-4t^2+7t+18}{12(t^2+t+1)}\cdot\beta_t & +\frac{-8t^2+10t}{12(t^2+t+1)}\cdot\beta_t^2 & +\frac{2t-3}{12(t^2+t+1)}\cdot\beta_t^3,\\
&&&\\
\frac{-4t^2+7t+18}{12(t^2+t+1)} & +\frac{40t^2+92t+62}{12(t^2+t+1)}\cdot\beta_t & +\frac{32t^2+46t+5}{12(t^2+t+1)}\cdot\beta_t^2 & +\frac{-8t-11}{12(t^2+t+1)}\cdot\beta_t^3,\\
&&&\\
\frac{-8t^2+10t}{12(t^2+t+1)} & +\frac{32t^2+46t+5}{12(t^2+t+1)}\cdot\beta_t & +\frac{32t^2+8t+1}{12(t^2+t+1)}\cdot\beta_t^2 & +\frac{-8t-1}{12(t^2+t+1)}\cdot\beta_t^3,\\
&&&\\
\frac{2t-3}{12(t^2+t+1)} & +\frac{-8t-11}{12(t^2+t+1)}\cdot\beta_t & +\frac{-8t-1}{12(t^2+t+1)}\cdot\beta_t^2 & +\frac{2}{12(t^2+t+1)}\cdot\beta_t^3
\end{array}
\right)$$

\begin{theorem}\label{periodl}
With the notation above, if $ n\not\equiv0\mod{3}$ and $t^2+t+1$ is square-free, else $t^2+3t+9$ is square-free, then an integral basis of the fields $\Q(\beta_t)$ is repeating periodically modulo $n_0$ given in the following table:
$$\renewcommand{\arraystretch}{1.5}
\scriptsize\begin{array}{c|c|c|c|c|c|c|c|c|c|c|c}
n & 2 & 3& 4& 5& 6& 7& 8& 9& 10& 11& 12\\\hline
n_0 & 2^2& 1^3& (3\cdot4)^4& (3^3\cdot5)^5& (3^2\cdot6)^6& (3^4\cdot7)^7& (3^5\cdot8)^8& (3^4\cdot9)^9& (3^6\cdot10)^{10}& (3^9\cdot11)^{11}& (3^8\cdot12)^{12}
\end{array}$$
\end{theorem}

\begin{prf}
In all cases, the least common multiple $d$ of the denominators of the elements of the dual basis represented in the $\Q$-basis $\left(1,\beta_t,\beta_t^2\ldots,\beta_t^{n-1}\right)$, has the form 
$$d=3^{\delta_n}\cdot n\cdot (t^2+t+1),\qquad \mbox{ if } n\not\equiv0\mod{3},$$
or 
$$d=3^{\delta_n}\cdot n \cdot (t^2+3t+9),\qquad \mbox{ if } n\equiv0\mod{3},$$
with
$$\begin{array}{c|c|c|c|c|c|c|c|c|c|c|c}
n & 2 & 3& 4& 5& 6& 7& 8& 9& 10& 11& 12\\\hline
\delta_n & 0& 0& 1& 3& 2& 4& 5& 4& 6& 9& 8
\end{array}$$
By Corollary \ref{denom} we can represent any algebraic integer in $\Q(\beta_t)$ in the $\Q$-basis $\left(1,\beta_t,\beta_t^2\ldots,\beta_t^{n-1}\right)$, having denominator $C_n$, for which,
$$C_n^2\mid 3^{\frac{n^2}{2}}\cdot n^n.$$
Thus, we can represent any algebraic integer $\gamma\in\Q(\beta_t)$ in the $\Q$-basis
$(1,\beta_t,\beta_t^2\ldots,\beta_t^{n-1})$, with rational coefficient having denominator 
$$\gcd(3^{\delta_n}\cdot n,C_n),$$
whence by Lemma \ref{cnn0}, an integral basis of these fields repeating periodically modulo $n_0$ given in the table.
\end{prf}

However these period lengths are much smaller than we give in Theorem \ref{period}, the smallest period lengths are even smaller. It was shown in \cite{sextic}, that an integral basis of the simplest sextic fields is repeating periodically modulo $36$ instead of $54^6$. In order to find the smallest period lengths, we have to calculate an integral basis of the fields belong to the parameters less than the period length given in the table above, and detect the $h_0^{(r)},\ldots, h_{n-1}^{(r)}$ polynomials, that are repeated. If the period length given in the table is too large for this calculations, we can use the method described in \cite{ibm} Theorem 3. For this we first have to make a guess for the smallest period length $n_0$, then we have to show, that there is no algebraic integer represented in the initial basis belonging to the residue class of our parameter, with rational coefficient having prime denominator $p$. We have to test it for all residue classes of $n_0$ and for all possible primes, which divides the initial period length given in the above table. For all of these primes this calculation requires to calculate an integral basis of  $n_0\cdot p^n$ number fields.\\
We can manage it for $n=2,3,4,5,6,8,9$ and $12$, but for $n=7,10,11$, the $7^7,5^{10},11^{11}$ tests would take too much time.\\
Therefore, we obtained the following result.

\begin{theorem}
With the notation above, if $ n\not\equiv0\mod{3}$ and $t^2+t+1$ is square-free, else $t^2+3t+9$ is square-free, then an integral basis of the fields $\Q(\beta_t)$ is repeating periodically modulo $n_0$ given in the following table:
$$\begin{array}{c|c|c|c|c|c|c|c|c}
n & 2 & 3& 4& 5& 6& 8& 9& 12\\\hline
n_0 & 4& 1& 24& 75& 36& 432& 1& 1944
\end{array}$$
\end{theorem}

\def\bibindent{4em}

\end{document}